\def\ung{{{\frak{g}}}}

\def\uqg{{{U_{q}(\hat{\ung}^\sigma)}}}
\def\uql{{{U_q(L(\ung)^\sigma)}}}

\def\unk{{\frak{k}}}
\def\bp{{{\bold{P}}}}
\def\bq{{{\bold{Q}}}}
\def\calp{{{{\Cal{P}}}}}

\def\ot{{{\otimes}}}

\def\calp{\Cal P}

\def\te0{\tilde{e}_0}
\def\xop{X_0^+}
\def\x1m{X_1^-}
\magnification 1200
\input amstex
\documentstyle{amsppt}
\NoBlackBoxes
\nologo
\document
\centerline{\bf{{Twisted Quantum Affine Algebras}}} 
\vskip36pt\centerline{Vyjayanthi Chari and Andrew Pressley}
\vskip 36pt
\noindent{\bf Introduction}  
\vskip12pt\noindent Quantum affine algebras are one of the most important classes of quantum groups. Their finite-dimensional representations lead to solutions of the quantum Yang--Baxter equation which are trigonometric functions of the spectral parameter (see [7], Sect. 12.5 B) and are thus related to various types of integrable models in statistical mechanics and field theory. Quantum affine algebras have also been shown to arise as \lq quantum symmetry groups\rq\ of certain integrable quantum field theories, such as affine Toda field theories (see [2] and [10]). More precisely, there is an affine Toda field theory associated to any affine Lie algebra $\unk$, and this theory admits as a quantum symmetry group the quantum affine algebra $U_q(\unk^*)$, where $\unk^*$ is the affine Lie algebra {\it dual} to $\unk$ (whose Dynkin diagram is obtained from that of $\unk$ by reversing the arrows). Since $\unk^*$ is often twisted even if $\unk$ is untwisted, this shows that the repr!
!
!
esentation theory of twisted quantum affine algebras is, in this context at least, just as important as that of untwisted ones. However, there appear to be virtually no results in the literature on the twisted case. The only exceptions appear to be [12] and [14], which prove the existence of finite-dimensional irreducible representations of twisted quantum affine algebras $U_q(\unk)$ which are irreducible under certain subalgebras of the form $U_q(\frak m)$, where 
$\frak m$ is a finite-dimensional Lie subalgebra of $\unk$, and [15] and [17] which construct, by vertex operator methods, quantum analogues of the standard modules (which are, of course, infinite-dimensional). 

In [6] and [8], we gave a classification of the finite-dimensional irreducible representations of untwisted quantum affine algebras in terms of their highest weights, which are in one-to-one correspondence with $n$-tuples of polynomials in one variable with constant coefficient one ($n$ being the rank of the underlying finite-dimensional Lie algebra). The purpose of this paper is to extend this result to the twisted case. We find that the finite-dimensional irreducible representations of twisted quantum affine algebras are again parametrized by $n$-tuples of polynomials. But $n$ is now the rank of the fixed point subalgebra of the diagram automorphism, and the way in which such an $n$-tuple determines a highest weight is more complicated than in the untwisted case. 

In the analogous classical situation, we classified in [8] the finite-dimensional irreducible representations of the twisted affine Lie algebra $\hat{\ung}^\sigma$, associated to a diagram automorphism $\sigma$ of a finite-dimensional complex simple Lie algebra $\ung$, by using the canonical embedding of $\hat{\ung}^\sigma$ in the untwisted affine Lie algebra $\hat\ung$. Namely, we showed that every finite-dimensional irreducible representation of $\hat\ung$ decomposes under $\hat{\ung}^\sigma$ into a finite direct sum of irreducibles, and that every finite-dimensional irreducible representation of $\hat{\ung}^\sigma$ arises in this way. Together with the results of [6] and [7], this gave the desired classification. In the quantum case, Jing [16] has shown how to embed $U_q(\hat{\frak g}^\sigma)$ into $U_q(\hat{\frak g})$, but this embedding is not as simple as in the classical case and we have preferred to use a direct approach, following the method used for untwisted quantum!
!
!
 affine algebras in [6] and [8]. Since the proofs are similar to those for the untwisted case, we omit many of the details. 

\vskip36pt\noindent{\bf 1 Twisted quantum affine algebras}
\vskip 12pt\noindent Let $\ung$ be a finite-dimensional complex simple Lie algebra with Cartan matrix $A=(a_{ij})_{i,j\in I}$. Let $\sigma:I\to I$ be a bijection such that $a_{\sigma(i)\sigma(j)}=a_{ij}$ for all $i,j\in I$, and let $m$ be the order of $\sigma$; we assume that $m>1$ (thus, $m=2$ or $3$). We also denote by $\sigma$ the corresponding Lie algebra automorphism of $\ung$. 

Fix a primitive $m$th root of unity $\omega\in\Bbb C^\times$. For $r\in\Bbb Z/m\Bbb Z$, let $\frak g_r$ be the eigenspace of $\sigma$ on $\ung$ with eigenvalue $\omega^r$. Then,
$$\ung=\bigoplus_{r\in\Bbb Z/m\Bbb Z}\ung_r$$
is a $\Bbb Z/m\Bbb Z$-gradation of $\ung$ (see [18], Chapter 8). 

The fixed point set $\ung_0$ of $\sigma$ is a simple Lie algebra. The nodes of its Dynkin diagram are naturally indexed by $I_\sigma$, the set of $\sigma$-orbits on $I$. Moreover, $\ung_1$ is an irreducible representation of $\ung_0$. Let $\{\alpha_i\}_{i\in I_\sigma}$ be a set of simple roots of $\ung_0$, and let $\theta$ be the highest weight of $\ung_1$ as a representation of $\ung_0$. Let $\{n_i\}_{i\in I_\sigma}$ be the positive integers such that 
$$\theta=\sum_{i\in I_\sigma}n_i\alpha_i.$$

The twisted affine Lie algebra $\hat{\ung}^\sigma$ is the universal central extension (with one-dimensional centre) of the twisted loop algebra 
$$L(\ung)^\sigma=\{f\in\Bbb C[t,t^{-1}]\otimes\ung\,\mid\,f(\omega t)=\sigma(f(t))\},$$
where $t$ is an indeterminate. It is well known (see [18]) that $\hat{\ung}^\sigma$ is a symmetrizable Kac--Moody algebra whose Dynkin diagram has nodes indexed by
$$\hat{I}_\sigma=I_\sigma\amalg\{0\}.$$
(Note: the node labelled $0$ here is labelled $\epsilon$ in [18].) Let $A^\sigma=(a_{ij}^\sigma)_{i,j\in\hat{I}_\sigma}$ be the (generalized) Cartan matrix of $\hat{\ung}^\sigma$, and let $\{d_i\}_{i\in\hat{I}_\sigma}$ be the coprime positive integers such that the matrix $(d_ia_{ij}^\sigma)$ is symmetric. Setting $n_0=1$, we have
$$\sum_{i\in\hat{I}_\sigma}n_id_ia_{ij}^\sigma=0\ \ \ \text{for all $j\in\hat{I}_\sigma$}.\tag1$$

Since $\hat{\ung}^\sigma$ is a symmetrizable Kac--Moody algebra there is, according to Drinfel'd and Jimbo, a corresponding quantum group $U_q(\hat{\ung}^\sigma)$. Namely, let $q$ be a non-zero complex number, {\it assumed throughout this paper not to be a root of unity}. Let $q_i=q^{d_i}$ for $i\in\hat{I}_\sigma$. If $n\in\Bbb Z$, set
$$[n]_q =\frac{q^n -q^{-n}}{q -q^{-1}},$$
and for $n\ge r\ge 0$, 
$$\align [n]_q! &=[n]_q[n-1]_q\ldots [2]_q[1]_q,\\
\left[{n\atop r}\right]_q &= \frac{[n]_q!}{[r]_q![n-r]_q!}.\endalign$$

\proclaim{Proposition 1.1} There is a Hopf algebra $\uqg$ over $\Bbb C$ which is generated as an algebra by elements $e_i^{{}\pm{}}$, $k_i^{{}\pm 1}$ ($i\in \hat{I}_\sigma$), with the following defining relations:
$$\align
k_ik_i^{-1} = k_i^{-1}k_i &=1;\;\;  k_ik_j =k_jk_i;\\
k_ie_j^{{}\pm{}}k_i^{-1} &= q_i^{{}\pm a_{ij}^\sigma}e_j^{{}\pm};\\
[e_i^+ , e_j^-] &= \delta_{ij}\frac{k_i - k_i^{-1}}{q_i -q_i^{-1}};\\
\sum_{r=0}^{1-a_{ij}^\sigma}(-1)^r\left[{{1-a_{ij}^\sigma}\atop r}\right]_{q_i} (e_i^{{}\pm{}})^re_j^{{}\pm{}}&(e_i^{{}\pm{}})^{1-a_{ij}^\sigma-r} =0\ \ {\text{if}} \ \ \ \ i\ne j.\endalign$$

The comultiplication $\Delta$ of $\uqg$ is given by
$$\Delta(e_i^+)= e_i^+\ot k_i +1\ot e_i^+,\ \ 
\Delta(e_i^-)= e_i^-\ot 1 +k_i^{-1}\ot e_i^-,\ \ 
\Delta(k_i^{{}\pm 1}) = k_i^{{}\pm 1}\ot k_i^{{}\pm 1}.$$
\endproclaim

It follows from (1) that
$$c=\prod_{i\in\hat{I}_\sigma}k_i^{n_i}$$
lies in the centre of $\uqg$. Let $\uql$ be the quotient of $\uqg$ by the ideal generated by $c-1$. Note that, since $c$ is group-like, $\uql$ inherits a natural Hopf algebra structure from $\uqg$. 

The following theorem is an analogue of a result of Drinfel'd ([13], Theorem 4).

Let $u_1$ and $u_2$ be independent indeterminates and, for $i,j\in I$, define $p_{i},d_{ij}\in\Bbb Q$, $P_{ij}^\pm,F_{ij}^\pm,G_{ij}^\pm\in\Bbb C[u_1,u_2]$ as follows:
\vskip12pt
if $\sigma(i)=i$, then $p_i=m$, $d_{ij}=\frac12$, $P_{ij}^\pm(u_1,u_2)=1$;
\vskip6pt
if $a_{i\sigma(i)}=0$ and $\sigma(j)\ne j$, then $p_i=1$, $d_{ij}=\frac1{4m}$, $P_{ij}^\pm(u_1,u_2)=1$;
\vskip6pt
if $a_{i\sigma(i)}=0$ and $\sigma(j)=j$, then $p_i=1$, $d_{ij}=\frac12$, $P_{ij}^\pm(u_1,u_2)=\frac{u_1^mq^{\pm 2m}-u_2^m}{u_1q^{\pm 2}-u_2}$;
\vskip6pt
if $a_{i\sigma(i)}=-1$, then $p_i=1$, $d_{ij}=\frac18$, $P_{ij}^\pm(u_1,u_2)=u_1q^{\pm 1}+u_2$;
\vskip6pt\hskip1cm
$F_{ij}^\pm(u_1,u_2)=\prod_{r\in\Bbb Z/m\Bbb Z}(u_1-\omega^rq^{\pm a_{i\sigma^r(j)}}u_2)$;
\vskip6pt\hskip1cm
$G_{ij}^\pm(u_1,u_2)=\prod_{r\in\Bbb Z/m\Bbb Z}(u_1q^{\pm a_{i\sigma^r(j)}}-\omega^ru_2)$.

\proclaim{Definition 1.2} For $i\in I$, let $\overline{i}\in I_\sigma$ be the $\sigma$-orbit of $i$. Let $D_q(\ung)^\sigma$ be the associative algebra over $\Bbb C$ with generators $X_{i,k}^\pm$ ($i\in I$, $k\in\Bbb Z$), $H_{i,k}$ ($i\in I$, $k\in\Bbb Z\backslash\{0\}$), $K_i^{\pm 1}$ ($i\in I$), and the following defining relations:

$$\align
X_{\sigma(i),k}^\pm=\omega^k X_{i,k}^\pm; & \ \ \ H_{\sigma(i),k}=\omega^k H_{i,k};\ \ \ K_{\sigma(i)}^{\pm 1}=K_i^{\pm 1};\\
K_iK_i^{-1}=K_i^{-1}K_i=1; & \ \ \ K_iK_j=K_jK_i;\\
H_{i,k}H_{j,l}=H_{j,l}H_{i,k}; & \ \ \ K_iH_{j,l}=H_{j,l}K_i;\\
K_iX_{j,k}^\pm K_i^{-1}&=q^{\pm \frac{m}{p_ip_j}\sum_{r\in\Bbb Z/m\Bbb Z}a_{i\sigma^r(j)}}X_{j,k}^\pm;\\
[H_{i,k},X_{j,l}^\pm]&=\pm\frac1k 
\left(\sum_{r\in\Bbb Z/m\Bbb Z}[ka_{i\sigma^r(j)}/d_{\overline{i}}]_{q_{\overline{i}}}\,\omega^{kr}\right)
X_{j,k+l}^\pm;\\
[X_{i,k}^+,X_{j,l}^-]&=\sum_{r\in\Bbb Z/m\Bbb Z}
\delta_{\sigma^r(i),j}\omega^{rl}
\left(\frac{\Psi_{i,k+l}^+-\Psi_{i,k+l}^-}{q_{\overline{i}}
-q_{\overline{i}}^{-1}}\right),
\endalign$$
where the $\Psi_{i,k}^\pm$ are defined by
$$\sum_{k=0}^\infty \Psi_{i,\pm k}^\pm u^{k}=K_i^{\pm 1}\text{exp}
\left(\pm(q_{\overline{i}}-q_{\overline{i}}^{-1})\sum_{l=1}^\infty H_{i,\pm l}u^{l}\right),$$
$u$ being an indeterminate, and $\Psi_{i,k}^\pm=0$ if $\mp k>0$;
$$F_{ij}^\pm(u_1,u_2)X_i^\pm(u_1)X_j^\pm(u_2)=G_{ij}^\pm(u_1,u_2)X_j^\pm(u_2)X_i^\pm(u_1),$$
where
$$X_i^\pm(u)=\sum_{k=0}^\infty X_{i,k}^\pm u^{-k};$$
$$\align
\text{Sym}\{P_{ij}^\pm(u_1,u_2)(X_j^\pm(v)X_i^\pm(u_1)X_i^\pm(u_2) &-(q^{2md_{ij}}+q^{-2md_{ij}})X_i^\pm(u_1)X_j^\pm(v)X_i^\pm(u_2)\\
&+X_i^\pm(u_1)X_i^\pm(u_2)X_j^\pm(v))\}=0
\endalign$$
if $a_{ij}=-1$ and $\sigma(i)\ne j$, where $u_1$, $u_2$ and $v$ are independent indeterminates and Sym denotes symmetrization over $u_1,u_2$;
$$\text{Sym}\{(q^{3/2}u_1^{\mp 1}-(q^{1/2}+q^{-1/2})u_2^{\mp 1}+q^{-3/2}u_3^{\mp 1})X_i^\pm(u_1)X_i^\pm(u_2)X_i^\pm(u_3)\}=0\tag{$2^\pm$} $$
and
$$\text{Sym}\{(q^{-3/2}u_1^{\pm 1}-(q^{1/2}+q^{-1/2})u_2^{\pm 1}+q^{3/2}u_3^{\pm 1})X_i^\pm(u_1)X_i^\pm(u_2)X_i^\pm(u_3)\}=0\tag{$3^\pm$}$$
if $a_{i\sigma(i)}=-1$, where Sym denotes symmetrization over the independent indeterminates $u_1,u_2,u_3$.\endproclaim

\proclaim{Theorem 1.3} There exists an isomorphism of algebras between $\uql$ and $D_q(\ung)^\sigma$ such that 
$$e_{\overline{i}}^+=X_{i,0}^+,\ \ \ e_{\overline{i}}^-=\frac1{p_i}X_{i,0}^-,\ \ \ k_{\overline{i}}=K_i,$$
where $i\in I$ belongs to the $\sigma$-orbit $\overline{i}$.\endproclaim

\vskip6pt\noindent{\it Remarks} 1. There is a similar realization of $U_q(\hat{\ung}^\sigma)$. Theorem 1.3 is, however, sufficient for our purposes since it can be shown (cf. [7], Proposition 12.2.3) that the central element $c$ acts as one on every finite-dimensional representation of $U_q(\hat{\ung}^\sigma)$.

2. Relations (2) and (3) are present only when $\hat{\ung}^\sigma$ is of type $A_{2n}^{(2)}$. Drinfel'd ([13], Theorem 4) has analogues of only two of these four relations (namely ($2^-$) and ($3^+$)). The other two can be shown to be consequences of these together with the other defining relations of $D_q(\ung)^\sigma$. We have included all four partly for reasons of symmetry, and partly because they are all needed in subsequent calculations.

3. The isomorphism in 1.3 depends on the choice of the section $\bar i\mapsto i$ of the canonical projection $I\to I_\sigma$, but any two such isomorphisms differ only by a rescaling on the generators $e_{\bar i}^\pm$ ($\bar i\in\hat{I}_\sigma$).

\vskip12pt 
For a proof of this theorem, and an explicit description of the isomorphism, see [16], Theorem 3.1 and [17], Proposition 2.1. However, in the $A_{2n}^{(2)}$ case, the $q$ in [16] and [17] must be replaced by $q^2$ to get the algebras denoted here by $\uql$ and $D_q(\ung)^\sigma$. Compare also [1] and [11] for analogous results in the untwisted case.

For later use, we record here the defining relations, and the form of the isomorphism in 1.3, for the simplest twisted quantum affine algebra $U_q(L(sl_3)^\tau)$, where $\tau$ is the non-trivial diagram automorphism of $sl_3(\Bbb C)$. In this case, we may drop the index $i$ from the generators of $U_q(L(sl_3)^\tau)$ (since $|I_\tau|=1$). The generalized Cartan matrix is
$$A^\tau=\left(\matrix 2&-1\\-4&2\endmatrix\right),$$
so that $d_0=4$ and $d_1=1$. 

The defining relations are as follows:

$$\align
KK^{-1}=K^{-1}K=1,\ \ KH_k=H_kK,\ \ &H_kH_l=H_lH_k,\ \ KX_k^\pm K^{-1}=q^{\pm 2}X_k^\pm,\\
[X_k^+,X_l^-]=\frac{\psi_{k+l}^+-\psi_{k+l}^-}{q-q^{-1}},\ \ \text{where}\ \ 
\sum_{k=0}^\infty\psi_{\pm k}^\pm u^k&=K^{\pm 1}exp\left(\pm(q-q^{-1})\sum_{l=1}^\infty H_{\pm l}u^l\right),\\
[H_k,X_l^\pm]=\pm\frac{[2k]_q}{k}(q^{2k}+q^{-2k}&+(-1)^{k+1})X_{k+l}^\pm\ \ \ \text{if $k\ne 0$},\\
X_{k+2}^\pm X_l^\pm+(q^{\mp 2}-q^{\pm 4})X_{k+1}^\pm X_{l+1}^\pm -q^{\pm 2}X_k^\pm X_{l+2}^\pm &\\
=q^{\pm 2}X_l^\pm X_{k+2}^\pm &+(q^{\pm 4}-q^{\mp 2})X_{l+1}^\pm X_{k+1}^\pm-X_{l+2}^\pm X_k^\pm,\\
Sym(q^3X_{k\mp 1}^\pm X_l^\pm X_m^\pm -(q+q^{-1})X_k^\pm X_{l\mp 1}^\pm X_m^\pm &+q^{-3}X_k^\pm X_l^\pm X_{m\mp 1}^\pm )=0,\\
Sym(q^{-3}X_{k\pm 1}^\pm X_l^\pm X_m^\pm -(q+q^{-1})X_k^\pm X_{l\pm 1}^\pm X_m^\pm &+q^{3}X_k^\pm X_l^\pm X_{m\pm 1}^\pm )=0,\endalign$$
where $Sym$ means the sum over all permutations of $k,l,m$. 

The isomorphism in 1.3 is given by
$$\align
e_0^+=K^{-2}(X_0^-X_1^- -q^2X_1^-X_0^-),\ \ & 
e_0^-=\frac1{[4]_q^2}(X_{-1}^+X_0^+-q^{-2}X_0^+X_{-1}^+)K^2,\\
k_0=K^{-2},\ \ e_1^+=X_0^+,\ \ & e_1^-=X_0^-,\ \ k_1=K.
\endalign$$

\vskip12pt Let $U_+^\sigma$ (resp. $U_-^\sigma$, $U_0^\sigma$) be the subalgebras of $U^\sigma$ generated by the $X_{i,k}^+$ (resp. by the $X_{i,k}^-$, by the $\Psi_{i,k}^\pm$) for $i\in I$, $k\in\Bbb Z$. 

\proclaim{Proposition 1.4} $U^\sigma=U_-^\sigma.U_0^\sigma.U_+^\sigma$. 
\endproclaim 

The proof is straightforward.

\vskip 36pt\noindent
{\bf 2 Some subalgebras of $\uql$}
\vskip12pt\noindent The study of untwisted quantum affine algebras can be reduced, to some extent at least, to the case of quantum affine $sl_2$, by noting that any algebra of the former type can be generated by finitely many copies of the latter (see [1], Proposition 3.8). In the twisted case, one needs $U_q(L(sl_3)^\tau)$ in addition, where $\tau$ is the unique non-trivial diagram automorphism of $sl_3(\Bbb C)$.

We recall the definition of quantum affine $sl_2$:

\proclaim{Definition 2.1} $U_q(L(sl_2))$ is the associative algebra with generators $X_k^\pm$ ($k\in\Bbb Z$), $H_k$ ($k\in \Bbb Z\backslash\{0\}$), $K^{\pm 1}$, and the following defining relations:
$$\align
KK^{-1}=K^{-1}K=1;\ \ & KH_k=H_kK; \ \ H_kH_l=H_lH_k;\\
KX_k^\pm K^{-1}&=q^{\pm 2}X_k^\pm;\\
[H_k,X_l^\pm]&=\pm\frac1k [2k]_qX_{k+l}^\pm;\\
X_{k+1}^\pm X_l^\pm-q^{\pm 2}X_l^\pm X_{k+1}^\pm &=
q^{\pm 2}X_k^\pm X_{l+1}^\pm-X_{l+1}^\pm X_k^\pm;\\
[X_k^+,X_l^-]&=\frac{\Psi_{k+l}^+-\Psi_{k+l}^-}{q-q^{-1}},\endalign$$
where
$$\sum_{k=0}^\infty\Psi_{\pm k}^\pm u^k=K^{\pm 1}\text{exp}\left(\pm(q-q^{-1})\sum_{l=1}^\infty H_{\pm l}u^l\right),$$
and $\Psi_k^\pm=0$ if $\mp k>0$. \endproclaim

See [7], Theorem 12.2.1 -- we have set the central element ${\Cal C}^{1/2}$ equal to one.

The result we need is the following:

\proclaim{Proposition 2.2} Let $i\in I$.

(i) If $\sigma(i)\ne i$ and $a_{i\sigma(i)}\ne 0$, there is a homomorphism of algebras $\varphi_i:U_q(L(sl_3)^\tau)\to\uql$ such that
$$\varphi_i(X_k^\pm)=X_{i,k}^\pm,\ \ \varphi_i(H_k)=H_{i,k},\ \ \varphi_i(\Psi_k^\pm)=\Psi_{i,k}^\pm,\ \ \varphi_i(K)=K_i.$$

(ii) If $\sigma(i)\ne i$ and $a_{i\sigma(i)}=0$, there is a homomorphism of algebras $\varphi_i:U_q(L(sl_2))\to\uql$ such that
$$\varphi_i(X_k^\pm)=X_{i,k}^\pm,\ \ \varphi_i(H_k)=H_{i,k},\ \ \varphi_i(\Psi_k^\pm)=\Psi_{i,k}^\pm,\ \ \varphi_i(K)=K_i.$$

(iii) If $\sigma(i)=i$, there is a homomorphism of algebras $\varphi_i:U_{q^m}(L(sl_2))\to\uql$ such that
$$\align
\varphi_i(X_k^+)=\frac1m X_{i,mk}^+,\ \ \varphi_i(X_k^-)&=\frac{q_{\overline{i}}-q_{\overline{i}}^{-1}}
{q^m-q^{-m}}X_{i,mk}^-,\\
\varphi_i(H_k)=\frac{q_{\overline{i}}-q_{\overline{i}}^{-1}}
{q^m-q^{-m}}H_{i,mk},\ \ \varphi_i(\Psi_k^\pm)&=\Psi_{i,mk}^\pm,\ \ \varphi_i(K)=K_i.
\endalign$$
\endproclaim
\demo{Proof} Straightforward verification, using 1.3 and 2.1.\qed\enddemo
\vskip12pt\noindent{\it Remarks} 1. We have dropped the subscript $i$ from the generators of $U_q(L(sl_3)^\tau)$ in (i), since $|I_\tau|=1$.

2. In (iii), the generators $X_{i,k}^\pm$, $H_{i,k}$, $\Psi_{i,k}^\pm$ of $\uql$ vanish if $k$ is not a multiple of $m$.

3. We expect that the homomorphisms $\varphi_i$ are injective, but we shall not need this. 

\vskip36pt\noindent{\bf 3 Finite-dimensional representations} 
\vskip12pt\noindent A representation $V$ of $U^\sigma$ (i.e. a left $U^\sigma$-module) is said to be of type I if each $k_i$ ($i\in\hat{I}_\sigma$) acts semisimply on $V$ with eigenvalues which are integer powers of $q_i$. It is not difficult to show that every finite-dimensional irreducible representation of $U^\sigma$ can be obtained from a type I representation by twisting with an automorphism of $U^\sigma$ of the form $e_i^+\mapsto \epsilon_ie_i^+$, $e_i^-\mapsto e_i^-$, $k_i\mapsto \epsilon_i k_i$, where each $\epsilon_i=\pm 1$ (cf. [7], Proposition 12.2.3). 

If $V$ is a type I representation of $U^\sigma$, a vector $v\in V$ is said to be a highest weight vector if $v$ is annihilated by the $X_{i,k}^+$ for all $i\in I$, $k\in\Bbb Z$, and is a simultaneous eigenvector for the elements of $U_0^\sigma$. If, in addition, $V=U^\sigma.v$, then $V$ is said to be a highest weight representation. Moreover, if $\{\psi_{i,k}^\pm\}_{i\in I,k\in\Bbb Z}$ are the scalars such that
$$\Psi_{i,k}^\pm.v=\psi_{i,k}^\pm v,$$
the pair of $(I\times\Bbb Z)$-tuples $\pmb{\psi}^\pm=\{\psi_{i,k}^\pm\}_{i\in I,k\in\Bbb Z}$ is called the highest weight of $V$ (or the weight of $v$). Note that we necessarily have
$$\aligned
\psi_{i,k}^\pm &=0\ \ \text{if $\mp k>0$},\ \ \ \ \psi_{i,0}^+\psi_{i,0}^-=1,\\
\psi_{\sigma(i),k}^\pm &=\omega^k\psi_{i,k}^\pm\ \ \ \ \text{for all $i\in I$, $k\in\Bbb Z$}.
\endaligned\tag4$$
Conversely, by the usual Verma module construction, it is easy to show that, for any $\pmb{\psi}^\pm=\{\psi_{i,k}^\pm\}_{i\in I,k\in\Bbb Z}$ satisfying (4), 
there is, up to isomorphism, exactly one irreducible representation $V(\pmb{\psi}^\pm)$ with highest weight $\pmb{\psi}^\pm$. 

The following theorem is the main result of this paper:

\proclaim{Theorem 3.1} (i) Every finite-dimensional irreducible type I representation of $\uql$ is highest weight.

(ii) If $\pmb{\psi}^\pm=\{\psi_{i,k}^\pm\}_{i\in I,k\in\Bbb Z}$, the highest weight representation $V(\pmb{\psi}^\pm)$ of $\uql$ is finite-dimensional if and only if there exist polynomials $P_i\in\Bbb C[u]$ ($i\in I$) with constant coefficient one such that
$$\sum_{k=0}^\infty \psi_{i,k}^+u^k=\sum_{k=0}^\infty \psi_{i,-k}^-u^{-k}=
\cases
q^{m\text{deg}P_i}\frac{P_i(q^{-2m}u)}{P_i(u)}\ &\text{if $\sigma(i)\ne i$ and $a_{i\sigma(i)}\ne 0$,}\\
q^{\text{deg}P_i}\frac{P_i(q^{-2}u)}{P_i(u)}\ &\text{if $\sigma(i)\ne i$ and $a_{i\sigma(i)}= 0$,}\\
q^{m\text{deg}P_i}\frac{P_i(q^{-2m}u^m)}{P_i(u^m)}\ &\text{if $\sigma(i)=i$,}\endcases$$
in the sense that the first two terms are the Laurent expansions of the third term about $u=0$ and $u=\infty$, respectively. \endproclaim

The proof of (i) is straightforward (cf. [7], Proposition 12.2.3). The proof of (ii) will occupy the next two sections.
\vskip12pt\noindent{\it Remark} \,\,Since $\Psi_{\sigma(i),k}^\pm=\omega^k \Psi_{i,k}^\pm$, the polynomials $P_i$, if they exist, necessarily satisfy the condition
$$P_{\sigma(i)}(u)=P_i(\omega u).\tag5$$
\vskip12pt
Let $\Pi$ be the set of $I$-tuples of polynomials $P_i\in\Bbb C[u]$ with constant coefficient one satisfying (5). If $\bold P=\{P_i\}_{i\in I}\in\Pi$, we denote by $V(\bold P)$ the irreducible highest weight representation $V(\pmb{\psi}^\pm)$ of $U^\sigma$ (abusing notation), the relation between $\bold P$ and $\pmb{\psi}^\pm$ being as in 3.1. 

\proclaim{Proposition 3.2} Let $\bold P=\{P_i\}_{i\in I},\bold Q=\{Q_i\}_{i\in I}\in\Pi$, and let $v_{\bold P}\in V(\bold P)$, $v_{\bold Q}\in V(\bold Q)$ be highest weight vectors. Then, $v_{\bold P}\otimes v_{\bold Q}$ is a highest weight vector in $V(\bold P)\otimes V(\bold Q)$ of weight $\pmb{\phi}^\pm$, where $\pmb{\phi}^\pm$ is related to the $I$-tuple $\{P_iQ_i\}_{i\in I}$ in the same way as $\pmb{\psi}^\pm$ is related to $\{P_i\}_{i\in I}$ in 3.1. \endproclaim

This will be proved in Sections 4 and 5. The following corollary is immediate:

\proclaim{Corollary 3.3} Let the notation be as in 3.2, and denote the $I$-tuple $\{P_iQ_i\}_{i\in I}$ by $\bold P\otimes\bold Q$. Then, $V(\bold P\otimes\bold Q)$ is isomorphic as a representation of $\uql$ to a subquotient of $V(\bold P)\otimes V(\bold Q)$.\endproclaim

\vskip36pt\noindent{\bf 4 The $U_q(L(sl_3)^\tau)$ case}
\vskip12pt\noindent In this section, we prove 3.1 and 3.2 for  $U_q(L(sl_3)^\tau)$, where $\tau$ is the non-trivial diagram automorphism of $sl_3(\Bbb C)$, and we denote $U_q(L(sl_3)^\tau)$ by $U^\tau$. The explicit form of the generators and relations of $U^\tau$ was given at the end of Section 1. It will be convenient to set
$$\tilde{e}_0=X_0^-X_1^- -q^2X_1^-X_0^-,$$
so that, in the isomorphism in 1.3, $e_0^-$ is a scalar multiple of $\te0 K^2$, and to write
$$(X_k^\pm)^{(r)}=\frac{(X_k^\pm)^r}{[r]_q!},\ \ \ (\tilde{e}_0)^{(r)}=\frac{(\tilde{e}_0)^r}{[r]_{q^4}!}.$$
The crucial result for the proof of 3.1 in this case is the next proposition. 

\proclaim{Definition 4.1} Define elements $\{\calp_r\}_{r\in\Bbb N}$ in $U_0^\tau$ by $\calp_0=1$ and
$$\calp_r=-\frac1{1-q^{-4r}}\sum_{j=0}^{r-1}\Psi_{j+1}^+\calp_{r-j-1}K^{-1}.
\tag6$$
\endproclaim

If we introduce the formal power series
$$\calp(u)=\sum_{r=0}^\infty \calp_r u^r,\ \ \Psi^\pm(u)=\sum_{r=0}^\infty \Psi_{\pm r}^\pm u^{\pm r}$$
in an indeterminate $u$, then 4.1 is equivalent to saying that $\calp(u)$ has constant coefficient one and that
$$\Psi^+(u)=K\frac{\calp(q^{-4}u)}{\calp(u)}.$$

Let $X^+$ be the linear subspace of $U^\tau$ spanned by the $X_k^+$ for $k\in\Bbb Z$.
\proclaim{Proposition 4.2} For all $r\in\Bbb N$, we have the following congruences (mod $U^\tau X^+$):
\vskip6pt\noindent (i)${}_r$\ \ $(X_0^+)^{(2r+2)}(\tilde{e}_0)^{(r+1)}\equiv (-1)^{r+1}q^{-2(r+1)(2r+1)}[4]_q^{r+1}\calp_{r+1}K^{2r+2}$;
\vskip6pt\noindent(ii)${}_r$\ \ $(X_0^+)^{(2r+1)}(\tilde{e}_0)^{(r+1)}\equiv (-1)^{r}q^{-4r(r+1)}[4]_q^{r+1}\sum_{j=0}^r X_{j+1}^-\calp_{r-j}K^{2r+1}$;
\vskip6pt\noindent(iii)${}_r$\ \ $(X_0^+)^{(2r+1)}(\tilde{e}_0)^{(r+1)}\equiv q^{-6r+2}[4]_qKX_1^-(X_0^+)^{(2r)}(\tilde{e}_0)^{(r)}$
\vskip6pt\hskip6cm$
+q^{-8r+4}\frac{[4]_q}{[2]_q[3]_q}K^2[H_1,(X_0^+)^{(2r-1)}
(\tilde{e}_0)^{(r)}]$.
\endproclaim
\demo{Proof} All three congruences are easily checked when $r=0$. 

Assuming (iii)${}_r$, one deduces (ii)${}_{r}$ from (i)${}_{r-1}$, (ii)${}_{r-1}$ and (6), and then (i)${}_r$ follows from (ii)${}_r$ by multiplying on the left by $X_0^+$ and using (6) again.

{\vbox{
Thus, the main point is to prove (iii)${}_r$. For this, one needs identities (7)--(16) below: 
$$\align
&(X_0^+)^{(r)}  X_1^+=-q^{-3r}[r-1]_qX_1^+(X_0^+)^{(r)}+q^{-3r+3}X_0^+X_1^+
(X_0^+)^{(r-1)};\tag7\\
[H_1,(X_0^+)^{(r)}]&=[3]_q
\left\{\left(\frac{q^{-3r+3}+q^{-r+3}-q^{-r+1}-q^{-r-1}}{q-q^{-1}}\right)\right.
X_1^+(X_0^+)^{(r-1)}\\
& {\phantom{xxxxxxxxxxxxxxxx}}\left.+q^{-2r+4}X_0^+X_1^+(X_0^+)^{(r-2)}\right\};
\tag8\\
[(X_0^+)^{(r)},X_1^-]&=q^{-r+1}KH_1(X_0^+)^{(r-1)}\\
&{\phantom{xxxxx}}+
\left(\frac{q^{-2r+1}+q^{-2r-1}-q^{-2r+5}-q^{-4r+5}}{q-q^{-1}}\right)
KX_1^+(X_0^+)^{(r-2)}\\
&{\phantom{xxxxxxxxxxxxxx}}-q^{-3r+5}KX_0^+X_1^+(X_0^+)^{(r-3)};\tag9\\
[(X_0^+)^{(r)},\tilde{e}_0]&=q^{-r+3}[4]_qKX_1^-(X_0^+)^{(r-1)}
+q^{-2r+4}(q^2+q^{-2})K^2H_1(X_0^+)^{(r-2)}\\
&\ \ \ +q^{-3r+6}
\left(\frac{q^{-5}+q^{-3}-q^3-q^{-2r+3}}{q-q^{-1}}\right)K^2X_1^+
(X_0^+)^{(r-3)}\\
&\ \ \ \ \ \ \ \ -q^{-4r+8}K^2X_0^+X_1^+(X_0^+)^{(r-4)};\tag10\\
&\ \ \ \ \ \ \te0\x1m=q^4\x1m\te0\ ;\tag11\\
&[H_1,\te0^{(r)}]=-q^{-4r+5}(q-q^{-1})[3]_q[4]_q\te0^{(r-1)}(\x1m)^2\ ;\tag12\\
&[\xop,\te0^{(r)}]=q^{-4r+4}[4]_q\te0^{(r-1)}\x1m K\ ;\tag13\\
&[(\xop)^{(r+1)},\te0]=q^{-r+2}[4]_qK\x1m(\xop)^{(r)}+q^{-r}[2]_qK(\xop)^{(r)}\x1m\\
&\ \ \ \ \ \ \ +\frac{q^{-2r}}{[3]_q}K^2[H_1,(\xop)^{(r-1)}]+q^{-2r+3}(q-q^{-1})K^2H_1(\xop)^{(r-1)}\ ;\tag14\\
&\x1m\te0^{(r)}K\equiv\frac1{[4]_q}\xop\te0^{(r+1)}\ \ \ \text{(mod $U^\tau X^+$)}\ ;\tag15\\
&\te0^{(r-1)}(\x1m)^2\equiv\frac{q^{8r-2}}{[4]_q^2}(\xop)^2\te0^{(r+1)}K^{-2}
-\frac{q^{4r-2}}{[4]_q}\te0^{(r)}H_1\ \ \ \text{(mod $U^\tau X^+$)}\ .\tag16
\endalign$$
Identities (7)--(10) are proved successively by induction on $r$; (11) is a consequence of ($3^-$); (12) and (13) are proved by induction on $r$ using (11); (14) follows from (8) and (9); congruence (15) follows from (11) and (13); and (16) follows by a double application of (13).

Finally, to prove (iii)${}_r$, we compute
$$\align
[r+1]_{q^4}&(\xop)^{(2r+1)}\te0^{(r+1)}\\
&=\frac{q^{-2r+4}}{[2]_q}K\x1m(\xop)^{(4r)}\te0^{(r)}+q^{-2r}[2]_qK(\xop)^{(2r)}\x1m\te0^{(r)}\\
&\ \ \ \ +\frac{q^{-4r}}{[3]_q}K^2[H_1,(\xop)^{(2r-1)}]\te0^{(r)}+q^{-4r+3}(q-q^{-1})K^2H_1(\xop)^{(2r-1)}\te0^{(r)}\\
&\ \ \ \ \ \ \ \ \ \ \ \ \ \ \ \ \ \ \ \ \ \ \ \ \ \ \ \ \ \text{(by (14))}\endalign$$
}}
\vfill\eject
$$\align
&=\frac{q^{-2r+4}}{[2]_q}K\x1m(\xop)^{(2r)}\te0^{(r)}+q^{-2r}[2]_qK(\xop)^{(2r)}\x1m\te0^{(r)}\\
&\ \ \ \ +\frac{q^{-4r}}{[3]_q}K^2[H_1,(\xop)^{(2r-1)}\te0^{(r)}]-\frac{q^{-4r}}{[3]_q}K^2(\xop)^{(2r-1)}[H_1,\te0^{(r)}]\\
&\ \ \ +q^{-4r+3}(q-q^{-1})K^2[H_1,(\xop)^{(2r-1)}\te0^{(r)}]+q^{-4r+3}(q-q^{-1})K^2(\xop)^{(2r-1)}\te0^{(r)}H_1\\
&=\frac{q^{-2r+4}}{[2]_q}K\x1m(\xop)^{(2r)}\te0^{(r)}+q^{-2r}[2]_qK(\xop)^{(2r)}\x1m\te0^{(r)}\\
&\ \ \ \ +\frac{q^{-4r+6}}{[3]_q}K^2[H_1,(\xop)^{(2r-1)}\te0^{(r)}]-\frac{q^{-4r}}{[3]_q}K^2(\xop)^{(2r-1)}[H_1,\te0^{(r)}]\\
&\ \ \ \ +q^{-4r+3}(q-q^{-1})K^2(\xop)^{(2r-1)}\te0^{(r)}H_1\\
&\equiv \frac{q^{-2r+4}}{[2]_q}K\x1m(\xop)^{(2r)}\te0^{(r)}+
q^{-2r-2}\frac{[2]_q[2r+1]_q}{[4]_q}(\xop)^{(2r+1)}\te0^{(r+1)}\\
&\ \ \ +\frac{q^{-4r+6}}{[3]_q}K^2[H_1,(\xop)^{(2r-1)}\te0^{(r)}]-\frac{q^{-4r}}{[3]_q}K^2(\xop)^{(2r-1)}[H_1,\te0^{(r)}]\\
&\ \ \ +q^{-4r+3}(q-q^{-1})K^2(\xop)^{(2r-1)}\te0^{(r)}H_1\ \ \ \ \ \ \ \ \text{(mod $U^\tau X^+$)}\\
&\ \ \ \ \ \ \ \ \ \ \ \ \ \ \ \ \ \ \ \text{(by (15) applied to the second term)}\\
&\equiv \frac{q^{-2r+4}}{[2]_q}K\x1m(\xop)^{(2r)}\te0^{(r)}+
q^{-2r-2}\frac{[2]_q[2r+1]_q}{[4]_q}(\xop)^{(2r+1)}\te0^{(r+1)}\\
&\ \ \ +\frac{q^{-4r+6}}{[3]_q}K^2[H_1,(\xop)^{(2r-1)}\te0^{(r)}]
+q^{-8r+5}(q-q^{-1})[4]_qK^2(\xop)^{(2r-1)}\te0^{(r-1)}(\x1m)^2\\
&\ \ \  +q^{-4r+3}(q-q^{-1})K^2(\xop)^{(2r-1)}\te0^{(r)}H_1\ \ \ \ \ \ \ \ \text{(mod $U^\tau X^+$)}\\
&\ \ \ \ \ \ \ \ \ \ \ \ \ \ \ \ \ \ \ \text{(by (12) applied to the fourth term)}\\
&\equiv \frac{q^{-2r+4}}{[2]_q}K\x1m(\xop)^{(2r)}\te0^{(r)}+
q^{-2r-2}\frac{[2]_q[2r+1]_q}{[4]_q}(\xop)^{(2r+1)}\te0^{(r+1)}\\
&\ \ \ +\frac{q^{-4r+6}}{[3]_q}K^2[H_1,(\xop)^{(2r-1)}\te0^{(r)}]\\
&\ \ \ \ +q^{-8r+5}(q-q^{-1})[4]_qK^2(\xop)^{(2r-1)}\left(\frac{q^{8r-2}}{[4]_q^2}(\xop)^2\te0^{(r+1)}K^{-2}-\frac{q^{4r-2}}{[4]_q}\te0^{(r)}H_1\right)\\
&\ \ \ \ +q^{-4r+3}(q-q^{-1})K^2(\xop)^{(2r-1)}\te0^{(r)}H_1\ \ \ \ \ \ \ \ \text{(mod $U^\tau X^+$)}\\
&\ \ \ \ \ \ \ \ \ \ \ \ \ \ \ \ \ \ \ \text{(by (16) applied to the fourth term)}\\
&\equiv \frac{q^{-2r+4}}{[2]_q}K\x1m(\xop)^{(2r)}\te0^{(r)}+
q^{-2r-2}\frac{[2]_q[2r+1]_q}{[4]_q}(\xop)^{(2r+1)}\te0^{(r+1)}\\
&\ \ \ +\frac{q^{-4r+6}}{[3]_q}K^2[H_1,(\xop)^{(2r-1)}\te0^{(r)}]+
\frac{1-q^{-2}}{[4]_q}[2r+1]_q[2r]_q(\xop)^{(2r+1)}\te0^{(r+1)}\\
&\ \ \ \ \ \ \ \ \ \ \ \ \ \ \ \ \ \ \ \ \ \ \ \ \ \ \ \ \text{(mod $U^\tau X^+$)}.\endalign$$
Collecting the terms involving $(\xop)^{(2r+1)}\te0^{(r+1)}$ on the left-hand side and simplifying gives (iii)${}_r$.\qed\enddemo

Let $V$ be the finite-dimensional irreducible type I representation of $U^\tau$ with highest weight given by the pair of $\Bbb Z$-tuples $\{\psi_k^\pm\}_{k\in\Bbb Z}$, and let $v$ be a highest weight vector in $V$. We have
$$k_0.v=q^{4r_0}v,\ \ k_1.v=q^{r_1}v$$
for some $r_0,r_1\in\Bbb Z$. Note that $r_1=-2r_0$; in particular, $r_1$ is even. Write $r_0=-r$, $r_1=2r$ from now on.

Regarding $V$ as a representation of the $U_q(sl_2)$ subalgebra of $U^\tau$ generated by $e_1^\pm$ and $k_1^{\pm 1}$, $v$ is a highest weight vector (so that $r\ge 0$), and we have a direct sum decomposition
$$V\cong \bigoplus_{p\in\Bbb N}V_p^{n_p},$$
where $V_p$ is the irreducible representation of $U_q(sl_2)$ of dimension $p+1$ (and on which $k_1$ acts with eigenvalues in $q^{\Bbb Z}$), and the $n_p\ge 0$ are certain multiplicities. By 1.4, $n_p=0$ if $p$ is odd or $>2r$. Applying both sides of (i)${}_s$ in 4.2 to $v$, it follows that $\calp_s.v=0$ if $s>r$. Hence,
$$\calp(u).v=P(u)v,$$
where $P\in\Bbb C[u]$ is a polynomial with constant coefficient 1 and degree $\le r$. By the remarks following 4.1,
$$\Psi^+(u).v=q^{2r}\frac{P(q^{-4}u)}{P(u)}v.\tag17$$
Multiplying both sides of (ii)${}_r$ on the left by $X_{-n-1}^+$, where $n\in\Bbb N$, we see that
$$\sum_{j=n}^r\Psi_{j-n}^+\calp_{r-j}.v=\sum_{j=0}^n\Psi_{j-n}^-\calp_{r-j}.v\tag18$$
if $n\le r$, and 
$$\sum_{j=0}^r\Psi_{j-n}^-\calp_{r-j}.v=0\tag19$$
if $n>r$. By 4.1, (18) is equivalent to
$$q^{2r}q^{-4(r-n)}\calp_{r-n}.v=\sum_{j=0}^n\Psi_{j-n}^-\calp_{r-j}.v.\tag20$$
Equations (19) and (20) are together equivalent to
$$\Psi^-(u).v=q^{2r}\frac{P(q^{-4}u)}{P(u)}v.\tag21$$

By similar arguments, one shows the existence of a polynomial $Q\in\Bbb C[u]$ with constant coefficient one such that
$$\Psi^+(u).v=\Psi^-(u).v=q^{-2r}\frac{Q(q^4u^{-1})}{Q(u^{-1})}v.\tag22$$
Equating (21) and (22) and comparing the expansions in negative powers of $u$, we see that
$$q^{2r-4{\roman{deg}}P}=q^{-2r},$$
whence ${\roman{deg}}P=r$. Similarly, ${\roman{deg}}Q=r$ (in fact, it is clear that $Q(u)$ is a scalar multiple of $u^rP(u^{-1})$).

This completes the proof of the \lq only if\rq\ part of 3.1(ii) in the $U^\tau$ case. Before proving the \lq if\rq\ part, we prove 3.2 in the $U^\tau$ case. This depends on the following proposition.

\proclaim{Proposition 4.3} Let $k\ge 0$.
\vskip6pt\noindent (i) \ \ $\Delta(X_k^+)\equiv\sum_{j=0}^k X_{k-j}^+\otimes\Psi_j^+ +1\otimes X_k^+$\ \ (mod $U^\tau(X^+)^2\otimes U^\tau$);
\vskip6pt\noindent (ii) \ \ $\Delta(\Psi_k^+)\equiv\sum_{j=0}^k\Psi_j^+\otimes\Psi_{k-j}^+$\ \ (mod $U^\tau X^+\otimes U^\tau+U^\tau\otimes U^\tau X^+$).
\endproclaim
\demo{Proof} Making use of 1.1 and the isomorphism in 1.3, one computes that
$$\Delta(H_1)=H_1\otimes 1+1\ot H_1-(q-q^{-1})[2]_qX_0^+\ot X_1^-+q^{-1}(q-q^{-1})(X_0^+)^2\ot K^2e_0^+.$$
The formula in (i) now follows by an easy induction on $k$. Then (ii) follows from (i) by using
$$\Psi_k^+=(q-q^{-1})[X_k^+,X_0^-].\ \ \ \ \qed$$
\enddemo

Part (ii) of 4.3 implies that, when acting on a tensor product of two highest weight vectors, $\Psi^+(u)$ acts as a group-like element of the formal power series Hopf algebra $U^\tau[[u]]$. Proposition 3.2 follows.

To prove the \lq if\rq\ part of 3.1(ii) for $U^\tau$, it suffices by 3.3 to show that $V(P)$ is finite-dimensional when ${\roman{deg}}P=1$. This is accomplished in the next proposition.

\proclaim{Proposition 4.4} The following is a representation of $U^\tau$, for any $a\in\Bbb C^\times$:
$$X_k^+\mapsto a^k[2]_q\left(\matrix
0&1&0\\0&0&(-1)^kq^{2k}\\0&0&0\endmatrix\right),\ \ 
X_k^-\mapsto a^k\left(\matrix 0&0&0\\1&0&0\\0&(-1)^kq^{2k}&0\endmatrix\right),$$
$$H_k\mapsto a^k\frac{[2k]_q}{k}\left(\matrix q^{-2k}&0&0\\0&(-1)^k-q^{2k}&0\\0&0&(-1)^{k+1}q^{4k}\endmatrix\right),\ \ 
K\mapsto\left(\matrix q^2&0&0\\0&1&0\\0&0&q^{-2}\endmatrix\right),$$
$$\Psi_k^+\mapsto(q^2-q^{-2})a^k\left(\matrix 1&0&0\\0&(-1)^kq^{2k}-1&0\\0&0&(-1)^{k+1}q^{2k}\endmatrix\right)\ \ {\text{if $k>0$,}}$$
$$\Psi_k^-\mapsto -(q^2-q^{-2})a^k\left(\matrix 1&0&0\\0&(-1)^kq^{2k}-1&0\\0&0&(-1)^{k+1}q^{2k}\endmatrix\right)\ \ {\text{if $k<0$.}}$$
\endproclaim
\demo{Proof} Straightforward verification.\qed\enddemo

The representation defined in 4.4, say $V_a$, is clearly irreducible and of type I. Moreover, if $\{\psi_k^\pm\}_{k\in\Bbb Z}$ is its highest weight, we have
$$\sum_{k=0}^\infty \psi_k^+u^k=q^2+\sum_{k=1}^\infty(q^2-q^{-2})a^ku^k=q^2\frac{P(q^{-4}u)}{P(u)},
$$
where $P(u)=1-au$. Thus, we have exhibited a finite-dimensional irreducible type I representation of $U^\tau$ with highest weight given (as in 3.1) by an arbitrary polynomial of degree one. This completes the proof of 3.1(ii) in the $U^\tau$ case.

\vskip36pt\noindent
{\bf 5 The general case}
\vskip12pt\noindent In this section, we outline the proofs of 3.1 and 3.2 for an arbitrary twisted quantum affine algebra $U_q(L(\ung)^\sigma)$, which we denote by $U^\sigma$.

Suppose then that $V$ is a finite-dimensional highest weight representation of 
$U^\sigma$ with highest weight vector $v$ and highest weight $\{\psi_{i,k}^\pm\}_{i\in I,k\in\Bbb Z}$. We consider three cases, as in Proposition 2.2:
\vskip6pt\noindent {\it Case (i):} $\sigma(i)\ne i$ {\it and} $a_{i\sigma(i)}\ne 0$. Note that $m=2$ in this case. Using the homomorphism $\varphi_i$ described in 2.1(i), we can view $V$ as a representation of $U_q(L(sl_3)^\tau)$, and as such $v$ is still a highest weight vector in $V$. By the $U_q(L(sl_3)^\tau)$ case of 3.1, proved in the previous section, there exists $P_i\in\Pi$ such that
$$\sum_{k=0}^\infty \psi_{i,k}^+u^k=\sum_{k=0}^\infty \psi_{i,-k}^-u^{-k}=q^{2{\roman{deg}}P_i}\frac{P_i(q^{-4}u)}{P_i(u)}.$$
{\it Case (ii):} $\sigma(i)\ne i$ {\it and} $a_{i\sigma(i)}=0$. This time, we can view $V$ as a representation of $U_q(L(sl_2))$. By [6], Theorem 3.4, there exists $P_i\in\Pi$ such that
$$\sum_{k=0}^\infty \psi_{i,k}^+u^k=\sum_{k=0}^\infty \psi_{i,-k}^-u^{-k}=q^{{\roman{deg}}P_i}\frac{P_i(q^{-2}u)}{P_i(u)}.$$
{\it Case (iii):} $\sigma(i)=i$. Viewing $V$ as a representation of $U_{q^m}(L(sl_2))$, there exists $P_i\in\Pi$ such that 
$$\sum_{k=0}^\infty \psi_{i,mk}^+u^k=\sum_{k=0}^\infty \psi_{i,-mk}^-u^{-k}=q^{m{\roman{deg}}P_i}\frac{P_i(q^{-2m}u)}{P_i(u)}.$$
Noting that $\Psi_{i,k}^\pm=0$ unless $k$ is divisible by $m$, we find that
$$\sum_{k=0}^\infty \psi_{i,k}^+u^k=\sum_{k=0}^\infty \psi_{i,-k}^-u^{-k}=q^{m{\roman{deg}}P_i}\frac{P_i(q^{-2m}u^m)}{P_i(u^m)}.$$

This proves the \lq only if\rq\ part of 3.1(ii). The \lq if\rq\ part is proved by an argument similar to that used in the untwisted case in [8], Sect. 5. In that case, the crucial point was to establish the result for $U_q(L(sl_2))$. In the present case, we also need the result for $U_q(L(sl_3)^\tau)$, which was proved at the end of Sect. 4. We omit further details.

Finally, to prove 3.3 in the general case, one uses the methods of [9], Sect. 2. Let $U_i$ denote $U_q(L(sl_3)^\tau)$, $U_q(L(sl_2))$ or $U_{q^m}(L(sl_2))$ in cases (i), (ii) or (iii) of 2.2, respectively. If $V$ is an irreducible highest weight representation of $U^\sigma$ with highest weight vector $v$, denote the representation $\varphi_i(U_i).v$ of $U_i$ by $V_i$.

\proclaim{Lemma 5.1} With the above notation, $V_i$ is an irreducible representation of $U_i$ with highest weight vector $v$.\endproclaim

The proof is similar to that of Lemma 2.3 in [9]. In particular, for any $\bold P=\{P_i\}_{i\in I}\in\Pi$, $V(\bold P)_i\cong V(P_i)$, where $V(P_i)$ is the finite-dimensional irreducible representation of $U_i$ associated to the polynomial $P_i$ as in 3.1(ii) if $U_i=U_q(L(sl_3)^\tau)$, and as in Theorem 3.4 in [6] if $U_i=U_q(L(sl_2))$ or $U_{q^m}(L(sl_2))$.

If $V$ and $W$ are two irreducible highest weight representations of $U^\sigma$ with highest weight vectors $v$ and $w$, then, for any $i\in I$, $V_i\ot W_i$ is a representation of $U_i$ via $(\varphi_i\ot\varphi_i)\circ\Delta_i$. On the other hand, it is not difficult to show that the subspace $V_i\ot W_i$ of $V\ot W$ is preserved by $(\Delta\circ\varphi_i)(U_i)$, giving a second way of viewing $V_i\ot W_i$ as a representation of $U_i$. We denote these representations by $V_i\ot_iW_i$ and $V_i\ot W_i$, respectively.

\proclaim{Lemma 5.2} With the above notation, the identity map $V_i\ot_iW_i\to V_i\ot W_i$ is an isomorphism of representations of $U_i$.\endproclaim

The proof is similar to that of Proposition 2.2 in [9]. The necessary facts about the comultiplication of $U_q(L(\frak g)^\sigma)$ can be established by computations similar to those used to prove Theorem 2.2 in [17]. 

Now let $\bp=\{P_i\}_{i\in I}$, $\bq=\{Q_i\}_{i\in I}\in\Pi$. Then, by the $U_q(L(sl_3)^\tau)$ case of 3.2, proved in the previous section, and the analogous result for $U_q(L(sl_2))$ (Proposition 4.3 in [6]), we have the following isomorphisms of representations of $U_i$:
$$V(\bp)_i\ot_i V(\bq)_i\cong V(P_i)\ot_i V(Q_i)\cong V(P_iQ_i).$$
If $v_{\bp}$ and $v_{\bq}$ are highest weight vectors in $V(\bp)$ and $V(\bq)$, it follows from 5.2 that
$$\Psi_{i,k}^\pm.(v_{\bp}\ot v_{\bq})=\psi_{i,k}^\pm(v_{\bp}\ot v_{\bq}),$$
the action on the left being given by $\Delta$, where $\{\psi_{i,k}^\pm\}_{i\in I,k\in\Bbb Z}$ corresponds to the polynomial $P_iQ_i$ as in the $U_q(L(sl_3)^\tau)$ case of 3.1(ii) (proved in the previous section) or the analogous result for $U_q(L(sl_2))$ or $U_{q^m}(L(sl_2))$  (Theorem 3.4 in [6]). This proves 3.2 for $U^\sigma$.

\vskip36pt\noindent{\bf References}
\vskip12pt
\noindent 1. Beck, J., Braid group action and quantum affine algebras, Commun. Math. Phys. {\bf 165} (1994), 555--568.

\noindent 2. Bernard, D. and LeClair, A., Quantum group symmetries and non-local currents in 2D QFT, Commun. Math. Phys. {\bf 142} (1991), 99--138.

\noindent 3. Chari, V., Integrable representations of affine Lie algebras, Invent. Math. {\bf 85} (1986), 317--335.

\noindent 4. Chari, V. and Pressley, A. N., New unitary representations of loop groups, Math. Ann. {\bf 275} (1986), 87--104.

\noindent 5. Chari, V. and Pressley, A. N., Integrable representations of twisted affine Lie algebras, J. Algebra {\bf 113} (1988), 438--464.

\noindent 6. Chari, V. and Pressley, A. N., Quantum affine algebras, Commun. Math. Phys. {\bf 142} (1991), 261--283.

\noindent 7.  Chari, V. and Pressley, A. N., {\it A Guide to Quantum Groups}, Cambridge University Press, Cambridge, 1994.

\noindent 8.  Chari, V. and Pressley, A. N., Quantum affine algebras and their representations, Canadian Math. Soc. Conf. Proc. {\bf 16} (1995), 59--78.

\noindent 9. Chari, V. and Pressley, A. N., Minimal affinizations of representations of quantum groups: the simply-laced case, J. Algebra {\bf 184} (1996), 1--30.

\noindent 10. Chari, V. and Pressley, A. N., Yangians, integrable quantum systems and Dorey's rule, Commun. Math. Phys. {\bf 181} (1996), 265--302.

\noindent 11. Damiani, I., A basis of type Poincar\'e--Birkhoff--Witt for the quantum algebra of $\hat{sl}(2)$, J. Algebra {\bf 161} (1993), 291--310.

\noindent 12. Delius, G. W., Gould, M. D. and Zhang, Y.-Z., Twisted quantum affine algebras and solutions to the Yang--Baxter equation, preprint KCL-TH-95-8, q-alg/9508012.

\noindent 13. Drinfel'd, V. G., A new realization of Yangians and quantized affine algebras, Soviet Math. Dokl. {\bf 36} (1988), 212--216.

\noindent 14. Gandenberger, G. M., McKay, N. J. and Watts, G. M. T., Twisted algebra R-matrices and S-matrices for $b_n^{(1)}$ affine Toda solitons and their bound states, preprint DAMTP-95-49, CRM-2314, hep-th/9509007.

\noindent 15. Jing, N.-H., Twisted vertex representations of quantum affine algebras, Invent. Math. {\bf 102} (1990), 663--690.

\noindent 16. Jing, N.-H., On Drinfeld realization of quantum affine algebras, preprint q-alg/9610035.

\noindent 17. Jing, N.-H. and Misra, K. C., Vertex operators for twisted quantum affine algebras, preprint q-alg/9701034. 

\noindent 18. Kac, V. G., {\it Infinite dimensional Lie algebras}, 3rd edition, Cambridge University Press, Cambridge, 1990.

\vskip36pt
{\eightpoint{
$$\matrix\format\l&\l&\l&\l\\
\phantom{.} & {\text{Vyjayanthi Chari}}\phantom{xxxxxxxxxxxxx} & {\text{Andrew\ Pressley}}\\
\phantom{.}&{\text{Department of Mathematics}}\phantom{xxxxxxxxxxxxx} & {\text{Department of Mathematics}}\\
\phantom{.}&{\text{University of California}}\phantom{xxxxxxxxxxxxx} & {\text{King's College}}\\
\phantom{.}&{\text{Riverside}}\phantom{xxxxxxxxxxxxx} & {\text{Strand}}\\
\phantom{.}&{\text{CA 92521}}\phantom{xxxxxxxxxxxxx} & {\text{London WC2R 2LS}}\\
\phantom{.}&{\text{USA}}\phantom{xxxxxxxxxxxxx} & {\roman{UK}}\\
&{\text{email: chari\@ucrmath.ucr.edu}}\phantom{xxxxxxxxxxxxx} & {\text{email:anp\@mth.kcl.ac.uk}}
\endmatrix$$
}}

\vfill\eject

\enddocument